%% file: main.tex
\documentclass[11pt,letterpaper]{article}
\input{sections/preamble}

\usepackage[margin=1in]{geometry}

\makeatletter
\def\blfootnote{\gdef\@thefnmark{}\@footnotetext}
\makeatother

\title{Negative Momentum for Convex-Concave Optimization}

	 \author{
	 	Henry Shugart
	 	\\	UPenn \\	\texttt{hshugart@upenn.edu}
	 	\and
        Shuyi Wang\\
  Yale\\
  \texttt{shuyi.wang@yale.edu}\\ 
  \and
	 	Jason M. Altschuler
	 	\\	UPenn \\	\texttt{alts@upenn.edu}
	 }

\begin{document}

\maketitle

\input{sections/abstract}

\input{sections/introduction}

 \input{sections/preliminaries}

\input{sections/main_results}

\input{sections/proof_key_lemma}

\newpage
\appendix
\input{sections/appommitteddetails}

\small
\bibliographystyle{plainnat}
\bibliography{refs}
\normalsize

\end{document}

%% file: sections/preamble.tex
\usepackage{amsmath,amssymb,amsthm}
\numberwithin{equation}{section}
\usepackage{MnSymbol} 
\usepackage{bm, bbm} 
\usepackage{graphicx,color}
\usepackage[hyphens]{url}
\usepackage{dsfont}
\usepackage{booktabs, makecell, adjustbox}
\usepackage[normalem]{ulem}
\usepackage{mathtools}
\usepackage{thmtools}
\usepackage[hidelinks,colorlinks=true,linkcolor=blue,citecolor=blue]{hyperref}
\usepackage[nameinlink,capitalize]{cleveref}
\usepackage{multirow}
\usepackage{algorithm}
\usepackage[noend]{algpseudocode}
\usepackage{tikz} 
\usepackage{pifont}
\usepackage{physics}
\usepackage{enumitem}

\usepackage{enumitem}

\usepackage{soul} 
\usepackage[square,sort,numbers]{natbib}

\usepackage{subcaption}

\usepackage{tabularx}
\newcolumntype{L}{>{\raggedright\arraybackslash}X}

\newtheorem{theorem}{Theorem}[section]

\newtheorem{lemma}[theorem]{Lemma}
\newtheorem{remark}[theorem]{Remark}

\newtheorem{question}[theorem]{Question}

\newtheorem{defin}[theorem]{Definition}

\newcommand{\cO}{\mathcal{O}}

\newcommand{\R}{\mathbb{R}}

\newcommand{\N}{\mathbb{N}}

\DeclareMathOperator{\logeps}{\log\tfrac{1}{\eps}}

\newcommand*{\eps}{\varepsilon}

\renewcommand{\epsilon}{\varepsilon}

\crefname{section}{\mbox{\S\!\!}}{\mbox{\S\!\!}}

\newcommand{\Msmooth}{M_{\mathrm{smooth}}}
\newcommand{\Mconvex}{M_{\mathrm{convex}}}
\newcommand{\Mconcave}{M_{\mathrm{concave}}}
\newcommand{\Mmonotone}{M_{\mathrm{monotone}}}

%% file: sections/abstract.tex
\begin{abstract}
    This paper revisits momentum in the context of min-max optimization. Momentum is a celebrated mechanism for accelerating gradient dynamics in settings like convex minimization, but its direct use in min-max optimization makes gradient dynamics diverge. Surprisingly,~\citep{gidel2019negative} showed that \textit{negative} momentum can help fix convergence. However, despite these promising initial results and progress since, the power of momentum remains unclear for min-max optimization in two key ways. (1) Generality: is global convergence possible for the foundational setting of convex-concave optimization? This is the direct analog of convex minimization and is a standard testing ground for min-max algorithms. (2) Fast convergence: is accelerated convergence possible for strongly-convex-strong-concave optimization (the only non-linear setting where global convergence is known)? Recent work has even argued that this is impossible.
    We answer both these questions in the affirmative. Together, these results put negative momentum on more equal footing with competitor algorithms, and show that negative momentum enables convergence significantly faster and more generally than was known possible.

\end{abstract}

%% file: sections/introduction.tex
\section{Introduction}\label{sec:intro}
Efficiently solving min-max problems (also called saddle-point problems) of the form 
\begin{align}
\nonumber
    \min_x \max_y f(x,y)
\end{align}
is essential for high-impact applications across machine learning, game theory, robust optimization, distributed optimization, constrained optimization, and more~\citep{benzi2005numerical,von1947theory,Bertsekas_2014,boyd2011distributed,ben2009robust,hast2013pid,Madry_Makelov_Schmidt_Tsipras_Vladu_2019,Goodfellow_Pouget-Abadie_Mirza_Xu_Warde-Farley_Ozair_Courville_Bengio_2014}.
A major challenge for min-max optimization is that, despite close connections to standard optimization, many algorithmic phenomena are fundamentally different. An infamous example is that although gradient descent (GD) converges for convex optimization with appropriately chosen positive stepsizes, the direct analog is false for min-max optimization~\citep{korpelevich1976extragradient}.

This paper revisits \emph{momentum}, a powerful algorithmic technique in standard optimization whose potential remains somewhat unclear for min-max optimization. The traditional intuition behind momentum is to augment an iterative algorithm $z_{t+1} = h(z_t)$ to
$z_{t+1} = h(z_t) + \beta(z_{t} -z_{t-1})$ for some momentum parameter $\beta > 0$, in order to accelerate the dynamics along important directions.
Celebrated results dating back to Polyak~\citep{polyak1964some} and Nesterov~\citep{nesterov1983convex} establish this phenomenon formally for convex minimization, proving accelerated convergence rates in fundamental settings such as (strongly) convex, smooth objectives. In recent decades, momentum has been established more generally as a powerful mechanism for accelerating iterative algorithms across a broad swathe of optimization settings; see the survey~\citep{d2021acceleration}.

\paragraph*{Direct adaptation of momentum fails.} However, a fundamental challenge in min-max optimization is that the straightforward adaptation of momentum fails to make gradient dynamics converge. To explain this, first recall that the analog of gradient descent for min-max optimization is \emph{gradient-descent-ascent} (GDA), in which $x$ takes a descent step while $y$ takes an ascent step:
\begin{equation}\label{eq:gda}
\begin{aligned}
    x_{t+1} &= x_t - \eta \nabla_x f(x_t,y_t)\,, \\
    y_{t+1} &= y_t + \eta \nabla_y f(x_t,y_t)\,.
\end{aligned}
\end{equation}
The direct adaptation of momentum for GDA would then be 
\begin{equation}\label{eq:gda-momentum-simultaneous}
\begin{aligned}
    x_{t+1} &= x_t - \eta \nabla_x f(x_t,y_t) + \beta (x_t - x_{t-1}) \,, \\
    y_{t+1} &= y_t + \eta \nabla_y f(x_t,y_t) + \beta (y_t - y_{t-1})\,.
\end{aligned}
\end{equation}
Yet this algorithm~\eqref{eq:gda-momentum-simultaneous} fails to converge for any choice of momentum parameter $\beta > 0$, even in simple settings such as the $1$-dimensional unconstrained bilinear problem $\min_x \max_y \, xy$~\citep{gidel2019negative}.

\paragraph*{Negative momentum and alternation.} Intriguingly,
\emph{negative momentum} $\beta < 0$ can help remedy this issue---an influential idea proposed by the seminal paper~\citep{gidel2019negative}.
This is counterintuitive from the perspective of standard (non min-max) optimization: negative momentum downweights past motion rather than upweights.
In min-max optimization, however, negative momentum helps because it partially suppresses the non-convergent rotational dynamics of GDA, so that the iterates spend less effort cycling around the saddle and more effort moving in genuinely improving directions.
Crucial for making this negative momentum work, even in the simple setting of bilinear $f$, is to \emph{alternate} the $x$ and $y$ updates:
\begin{equation}\label{eq:gda-momentum-alternating}
\begin{aligned}
    x_{t+1} &= x_t - \eta \nabla_x f(x_t,y_t) + \beta (x_t - x_{t-1}) \,, \\
    y_{t+1} &= y_t + \eta \nabla_y f(x_{t+1},y_t) + \beta (y_t - y_{t-1}).
\end{aligned}
\end{equation}
(Note the use of $x_{t+1}$ to update $y_{t+1}$.)
Specifically,~\citep{gidel2019negative} showed theoretically that this algorithm~\eqref{eq:gda-momentum-alternating} (with negative momentum and alternation) enables global convergence for bilinear objectives $f$, and helps empirically for complicated non-convex-non-concave settings like GANs.

However, despite these initial promising results and progress in nearly the decade since, the power of negative momentum remains unclear for two key reasons: \emph{generality} and \emph{slow convergence}.

\paragraph*{Generality of negative momentum?} Indeed, it has proven quite difficult to establish global\footnote{Of course, local convergence rates can be shown by linearizing the dynamics in order to reduce to the setting where $\nabla f$ is linear~\citep{gidel2019negative}. However, such arguments are unable to give non-asymptotic global convergence rates.} convergence rates for negative momentum for settings beyond linear $\nabla f$ (i.e., beyond bilinear and quadratic objectives $f$) because then the dynamics of the algorithm are non-linear and more complicated. In the min-max literature, the canonical proving ground for algorithms is convex-concave objectives $f$. Competitor algorithms (e.g., GDA with extragradients, optimism, etc.) are all classically known to converge in this setting \citep{korpelevich1976extragradient,Popov_1980}. In sharp contrast, it remains unknown whether negative momentum can ensure global convergence for this foundational setting.

\begin{question}[Convergence for convex-concave $f$]\label{q:general}
    For convex-concave optimization, is any global convergence result true for GDA with negative momentum? 
\end{question}

Answering this question is essential for understanding the power of negative momentum, because there are many situations where optimization algorithms work well in linear settings but fail beyond. An infamous example is momentum in convex optimization. Indeed,~\citet{polyak1964some} showed in the 1960s that momentum accelerates GD on convex quadratics $f$ (i.e., linear $\nabla f$); however, it remained open for half a century whether momentum makes GD converge at a similarly accelerated rate for non-quadratic convex optimization (i.e., beyond linear $\nabla f$), and recently this was answered in the negative for Polyak's heavy ball method~\citep{Lessard_Recht_Packard_2016} and more generally for a wide range of momentum parameters~\citep{goujaud2025provable}. Does momentum similarly have fundamental failures for min-max optimization?

\paragraph*{Fast convergence of negative momentum?} Beyond linear $\nabla f$, convergence is known for negative momentum only in the setting of strongly-convex-strong-concave $f$~\citep{Zhang_Grosse}. This is a standard benchmark setting in the literature that serves as a stepping stone towards general convex-concave $f$. Competitor algorithms (e.g., GDA with extragradients, optimism, etc.) are all known to converge to an $\eps$-optimal solution in $\mathcal{O}(\kappa \logeps)$ iterations, where $\kappa$ denotes the condition number~\citep{mokhtari2020unified}. This rate is optimal among arbitrary first-order algorithms~\citep{Zhang_Hong_Zhang_2022}. In contrast, negative momentum is believed to be unable to achieve these fast rates: recent work~\citep{Zhang_Wang_2021} has argued that negative momentum is fundamentally suboptimal in that it requires at least $\Omega(\kappa^{1.5}\logeps)$ iterations, which is substantially slower.

\begin{question}[Optimal convergence for strongly-convex-strongly-concave $f$]\label{q:fast}
    For strongly-convex-strongly-concave optimization with condition number $\kappa$, can GDA with negative momentum converge faster than $\cO(\kappa^{1.5} \logeps)$? At the optimal rate of $\cO(\kappa \logeps)$?
\end{question}

A potential hope for circumventing the aforementioned lower bound $\Omega(\kappa^{1.5} \logeps)$ is that it applies only to negative momentum \emph{without} alternation, i.e., the algorithm~\eqref{eq:gda-momentum-simultaneous} rather than~\eqref{eq:gda-momentum-alternating}. As described above, alternation is necessary for negative momentum to converge for bilinear $f$~\citep{gidel2019negative}; however, it is not necessary for convergence under the strong growth conditions enjoyed by strongly-convex-strongly-concave $f$. 
A technical challenge for analyzing alternation is that it requires analyzing multi-step progress (at least $2$ updates rather than $1$), which is particularly difficult for non-linear $\nabla f$ (see also~\cref{q:fast}). This challenge has limited prior work on both upper and lower bounds, and therefore the power of negative momentum remains unclear.

\subsection{Contributions}\label{ssec:intro-cont}

This paper resolves both questions in the affirmative. Together, these results show that negative momentum enables GDA to converge substantially faster and more generally than was known possible. These results put negative momentum on more equal footing with competitor algorithms (e.g., extragradients, optimism, etc.) which are well-known to have positive answers to both questions.

\begin{theorem}[Informal version of~\cref{thm:cc-long}]\label{thm:cc-short}
    Suppose $f$ is convex-concave and smooth. With appropriate choices of the stepsize $\eta > 0$ and negative momentum $\beta < 0$, the algorithm~\eqref{eq:gda-momentum-alternating} converges in $\cO(\tfrac{1}{\eps})$ iterations. 
\end{theorem}

\begin{theorem}[Informal version of~\cref{thm:scsc-long}]\label{thm:scsc-short}
    Suppose $f$ is $\mu$-strongly-convex-strongly-concave and $L$-smooth. Let $\kappa := L/\mu$ denote the condition number. With appropriate choices of the stepsize $\eta > 0$ and negative momentum $\beta < 0$, the algorithm~\eqref{eq:gda-momentum-alternating} converges in $\cO(\kappa \logeps)$ iterations. 
\end{theorem}

\begin{figure}[t]
    \centering
    \begin{subfigure}[t]{0.49\linewidth}
        \centering
        \includegraphics[width=0.85\linewidth]{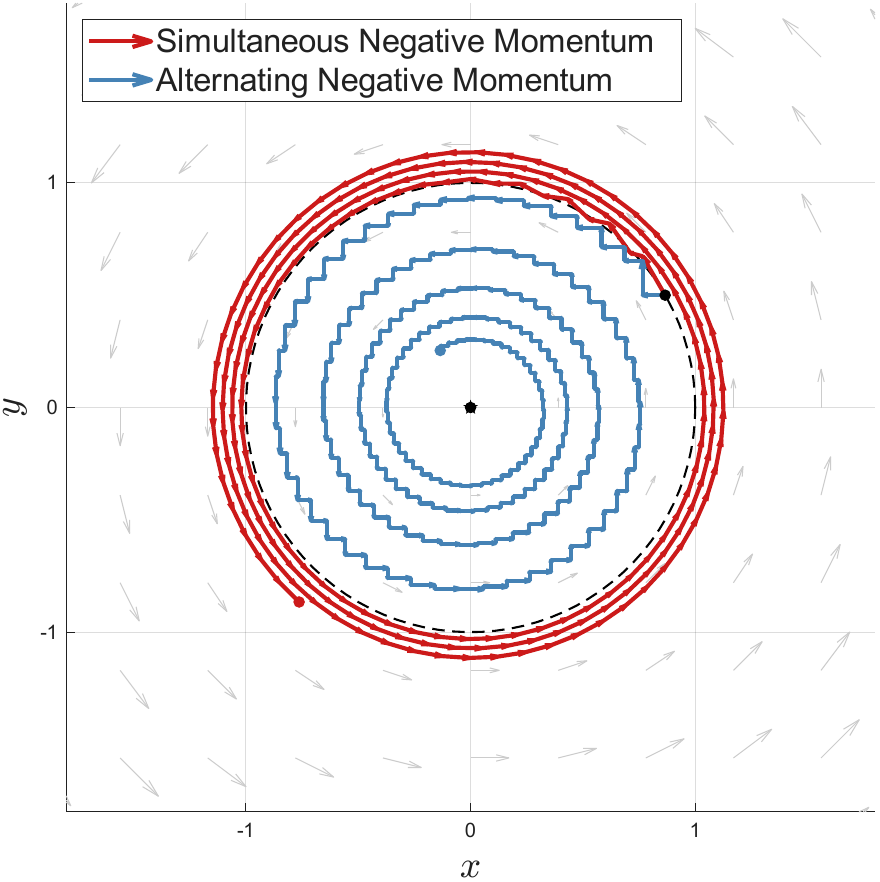}
        \caption{\footnotesize 
        For convex-concave $f$, alternation enables negative momentum to \emph{converge} (\cref{thm:cc-short}).
        }
    \end{subfigure}
    \hfill
    \begin{subfigure}[t]{0.49\linewidth}
        \centering
        \includegraphics[width=0.85\linewidth]{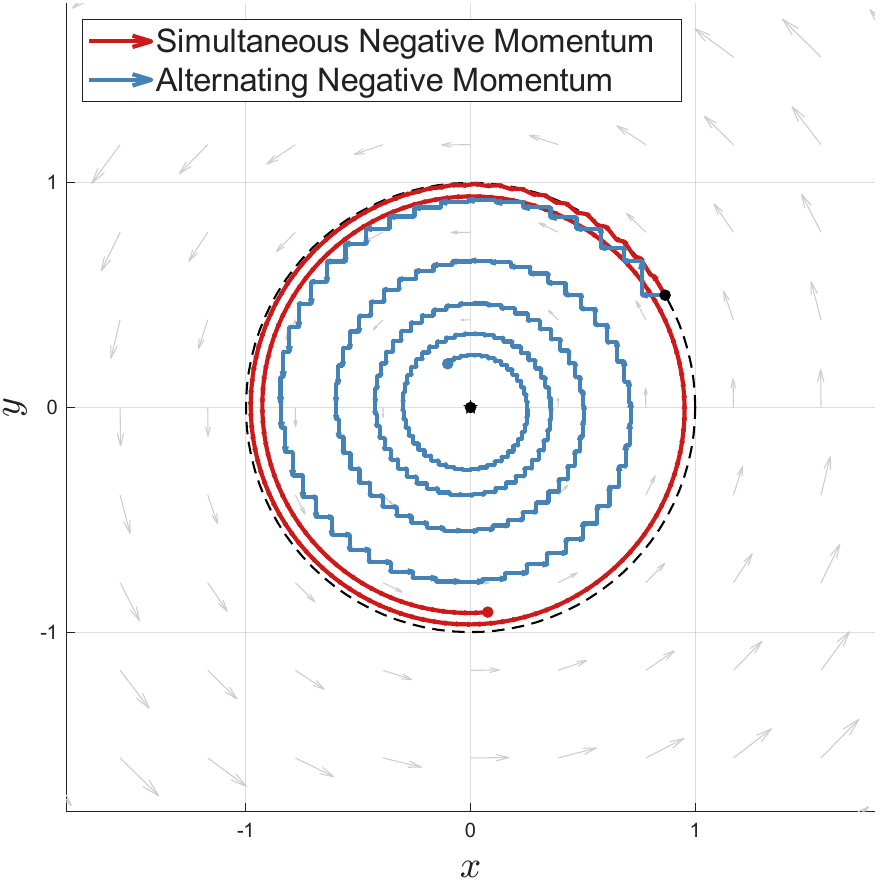}
        \caption{\footnotesize
        For strongly-convex-strongly-concave $f$, alternation enables negative momentum to \emph{accelerate} (\cref{thm:scsc-short}).
        }
    \end{subfigure}
    \caption{\footnotesize 
     Negative momentum is more effective with alternating updates~\eqref{eq:gda-momentum-alternating} in $x$ and $y$ than simultaneous updates~\eqref{eq:gda-momentum-simultaneous}. 
     % This makes the analysis more challenging (as it requires analyzing multi-step progress), but is essential for all our results.
     Plotted: 200 iterations on $1$-smooth problems $\min_x \max_y f(x,y)$ with unique solution at the origin. Alternating updates are shown for the parameter choices we analyze (stepsize $\eta = 0.2$ and momentum $\beta=-0.5$). Left: trajectories for convex-concave $f(x,y)=xy$. Simultaneous updates shown for $\eta=0.2$, $\beta=-0.8$; all parameter choices similarly diverge~\citep{gidel2019negative}. Right: trajectories for $\mu$-strongly-convex-strongly-concave $f(x,y)=\frac{\mu}{2} x^2 + \sqrt{1-\mu^2}xy-\frac{\mu}{2} y^2$, with $\mu = 0.01$. Simultaneous updates shown for the suggested optimal $\eta=\sqrt{\mu}=0.1, \beta=\sqrt{\mu}-1=-0.9$ from \citep{Zhang_Grosse}.
    }
    \label{fig:visual}
\end{figure}

Central to both results is the use of alternating updates in $x$ and $y$ rather than simultaneous updates (i.e., algorithm~\eqref{eq:gda-momentum-alternating} rather than~\eqref{eq:gda-momentum-simultaneous}). As described above, this makes the analysis more challenging as it requires  analyzing multi-step progress (at least $2$ updates rather than $1$). However, as we show, this alternation enables negative momentum to converge in convex-concave settings, and to accelerate in strongly-convex-strongly-concave settings (overcoming the aforementioned lower bound of~\citep{Zhang_Wang_2021} which applies only to simultaneous updates). See~\cref{fig:visual}.

We prove both results via a key progress lemma (\cref{lem:progress}) which identifies a non-obvious quadratic Lyapunov function under which negative momentum makes significant progress in each full iteration (i.e., after the update of both $x$ and $y$). This lemma applies to both settings.

\subsection{Related work}\label{ssec:intro-related}

\paragraph*{Negative momentum.} 
As described above, the direct use of (positive) momentum fails to make gradient dynamics converge for min-max optimization.~\citep{gidel2019negative} introduced negative momentum as an algorithmic mechanism for min-max problems, proving non-asymptotic global convergence for bilinear objectives and asymptotic local convergence for settings beyond, and showing empirically that this can help in complicated non-convex-non-concave settings like GANs. Their use of negative momentum was inspired by a trend in decreasing momentum values in GAN training.
In the years since, a number of works have studied negative momentum, e.g., for strongly-convex-strongly-concave objectives~\citep{Zhang_Wang_2021,Zhang_Grosse}, for constrained settings~\citep{fang2025rapid}, and
more broadly as an algorithmic building block~\citep{lorraine2022complex,feng2025continuous, beznosikov2022scaled}.
However, despite significant progress over the past decade, basic questions remain about the power of negative momentum, in particular about its generality (\cref{q:general}) and its speed (\cref{q:fast}). The purpose of this paper is to answer these two questions. 

\paragraph*{Asymmetry and alternation.} A key feature of some (but not all) min-max optimization algorithms is asymmetry, meaning updates that distinguish the minimization variable $x$ and maximization variable $y$. This asymmetry is provably beneficial in several contexts, 
for example accelerating the convergence of GDA in strongly-convex-strongly-concave settings~\citep{Lee_Cho_Yun_2024,Zhang_Wang_Lessard_Grosse_2022,shugart_negative_2025}, 
enabling the convergence of GDA in convex-concave settings~\citep{shugart_negative_2025}, and enabling faster algorithms for min-max problems than their associated variational inequalities~\citep{shugart2025min}. In the context of negative momentum, asymmetry is exploited through the alternation of updates of $x$ and $y$ (cf.,~\eqref{eq:gda-momentum-alternating} versus~\eqref{eq:gda-momentum-simultaneous}). This is provably necessary for convergence even for bilinear objectives~\citep{gidel2019negative} and therefore also for more general settings such as convex-concave objectives (cf.,~\cref{q:general}). For strongly-convex-strongly-concave objectives, alternation is necessary for negative momentum to converge at optimal rates~\citep{Zhang_Wang_2021}, and in this paper we show that it is also sufficient (cf.,~\cref{q:fast} and~\cref{thm:scsc-short}).

\paragraph*{Alternative algorithms.} Due to the failure of GDA~\eqref{eq:gda} even in simple bilinear settings such as $\min_x \max_y xy$~\citep{korpelevich1976extragradient}, an extensive literature has been devoted to developing convergent first-order algorithms for min-max optimization. Classic algorithms such as extragradient~\citep{korpelevich1976extragradient} and optimistic gradient~\citep{Popov_1980} enable GDA to converge for convex-concave objectives. 
These algorithms augment GDA by moving in the gradient direction at a ``lookahead iterate'' rather than at the current iterate.
While some qualitative connections have been made between optimism and negative momentum \citep{daskalakis2018limit, mokhtari2020unified}, connections thus far appear to be primarily syntactic and in particular no formal reductions or relations between theoretical convergence rates are known. 
Non-asymptotic convergence rates have recently been shown for these optimistic and extragradient algorithms, namely $\mathcal{O}(1/\eps)$ for convex-concave objectives \citep{Mokhtari_Ozdaglar_Pattathil_2020,cai2022tight,Gorbunov_Taylor_Gidel,gorbunov2022extragradient} and $\mathcal{O}(\kappa \log 1/\eps)$ for strongly convex-strongly concave objectives \citep{TSENG1995237, Daskalakis_Ilyas_Syrgkanis_Zeng_2018}. In this paper we show matching convergence rates for negative momentum, thereby putting it on more equal footing with these classical algorithms. For the strongly-convex-strongly-concave setting, this rate is optimal among arbitrary first-order algorithms. For the convex-concave setting, more sophisticated algorithms have recently been shown to accelerate convergence to $\mathcal{O}(1/\sqrt{\eps})$~\citep{Lee_Kim_2021, yoon2021accelerated, Tran_Dinh_Luo_2021}, and it remains an interesting problem if negative momentum can further accelerate on convex-concave objectives. By providing the first global convergence result for negative momentum,~\cref{thm:cc-long} opens the door to such refined questions.

\paragraph*{Computer-assisted search for Lyapunov functions.} Our convergence analysis is inspired by the framework of performance estimation problem (PEP), although with several key differences (detailed in~\cref{sec:proof}). Introduced in \citep{drori2014performance}, PEP casts the problem of finding the worst-case rate of a first-order algorithm over a fixed number of iterations as a semidefinite program. Closest to our work is the variant of PEP which automatically searches for Lyapunov functions certifying descent~\citep{taylor2018lyapunov, upadhyaya2025automated,Lessard_Recht_Packard_2016}. 
These families of ideas enable (partially) automating the design and analysis of optimization algorithms, and have been used to great effect in the past ${\sim}15$ years in a broad range of optimization settings (see the surveys~\citep{taylor2024towards, d2021acceleration} and references within), including min-max optimization (see for example \citep{Zhang_Wang_Lessard_Grosse_2022, krivchenko2026strengthening, krivchenko2024convex, Lee_Cho_Yun_2024, shugart_negative_2025, Ryu_Taylor_Bergeling_Giselsson_2020, Gorbunov_Taylor_Gidel, gorbunov2022extragradient, zamani2024convergence, yoon2021accelerated}).
We describe how our analysis framework builds upon PEP in detail in~\cref{sec:proof}.

%% file: sections/preliminaries.tex
\section{Preliminaries and notation}\label{sec:prelim}

We focus on unconstrained convex-concave optimization, a fundamental setting for min-max optimization that 
is a standard testing ground for algorithms. These problems are of the form
\begin{align}\label{eq:min-maxcc}
    \min_{x \in \R^{d_x}} \max_{y \in \R^{d_y}}\, f(x,y)\,.
\end{align}
Throughout we assume the existence of a solution, i.e., a point $(x^*,y^*)$ satisfying the stationarity condition $\nabla f(x^*,y^*) = 0$, or equivalently (since $f$ is convex-concave) satisfying the saddle-point condition $f(x^*,y) \leq f(x^*,y^*) \leq f(x,y^*)$ for all $x,y$. 
For simplicity of exposition, throughout we consider finite-dimensional spaces of (arbitrary) dimension $d_x$ and $d_y$, although we remark that our results extend to infinite-dimensional Hilbert space $\ell_2$.

Algorithmic results such as fast global convergence are tied to structural assumptions of the objective $f$. We focus on the standard setting of objectives $f$ that are smooth and convex-concave (or strongly-convex-strongly-concave). We briefly recall the definitions of these notions.

\begin{defin}[Smooth functions]\label{def:smooth}
    A function $g$ is $L$-smooth if its gradient $\nabla g$ is $L$-Lipschitz. That is, $\|\nabla g(a) - \nabla g(b)\| \leq L\|a-b\|$ for all $a,b$.
\end{defin}

\begin{defin}[(Strongly) convex functions]
A function $g$ is convex if $g(ta + (1-t)b) \leq tg(a) + (1-t)g(b)$ for all $t \in [0,1]$ and $a,b$. For $\mu \geq 0$, $g$ is $\mu$-strongly convex if $g(b)-\tfrac{\mu}{2}\|b\|^2$ is convex. 
\end{defin}

\begin{defin}[(Strongly) convex-concave functions]
For $\mu \geq 0$, a function $f(x,y)$ is $\mu$-strongly-convex-strongly-concave if $f(\cdot,y)$ is $\mu$-strongly-convex for every $y$, and $f(x,\cdot)$ is $\mu$-strongly-concave for every $x$. If $f$ satisfies this for $\mu=0$, $f$ is convex-concave.
\end{defin}

Our results do not require the objective $f \in C^2$. However, for intuition, we note that under such an assumption, $f$ being $L$-smooth is equivalent to $\|\nabla^2 f\|_{\mathrm{op}} \leq L$, and $f$ being $\mu$-strongly-convex-strongly-concave is equivalent to $\nabla_{xx}^2 f \succeq \mu \bm I_{d_x}$ and $-\nabla_{yy}^2 f \succeq \mu \bm I_{d_y}$.

In our analysis, the primary way that we exploit these structural properties of smoothness and (strong) convexity-concavity is through the standard ``co-coercivity inequality'', recalled next.

\begin{lemma}[Co-coercivity]\label{lem:cocoercivity}
Suppose $g$ is 1-smooth and $\mu$-strongly convex. The co-coercivity
\begin{align*}
    C_g(a,b) := g(a) - g(b) - \nabla g(b)^\top  (a-b) - \frac{\mu}{2}\|a-b\|^2 - \frac{1}{2(1-\mu)}\|\nabla g(a) - \nabla g(b) - \mu(a-b)\|^2
\end{align*}
satisfies $C_g(a,b) \geq 0$ for any points $a,b$.
\end{lemma}

\paragraph*{Notation.} Our notation is standard. We write $\bm I_d$ to denote the identity matrix in dimension $d$, $\otimes$ to denote the Kronecker product, and $\succeq,\preceq$ to denote inequalities in the L\"oewner order. The notation $\mathcal{O}(\cdot)$ and $\Omega(\cdot)$ refer to upper and lower bounds, respectively, up to universal constants. For shorthand, we often concatenate variables as $z_t := (x_t,y_t)$. 
All vectors are column vectors.

%% file: sections/main_results.tex
\section{Formal statement of main results}\label{sec:main}
In this section, we formally state our convergence results for alternating GDA with negative momentum on convex-concave and strongly-convex-strongly-concave functions.

Our first result answers~\cref{q:general}:

\begin{theorem}[Convergence for convex-concave objectives]\label{thm:cc-long}
    Let $f$ be an $L$-smooth, convex-concave function with saddle point $z^* = (x^*,y^*)$. Then for any dimensions $d_x,d_y \in \N$, any initialization $z_0 = (x_0,y_0)$, and any number of iterations $T \in \N$, the iterates of alternating GDA~\eqref{eq:gda-momentum-alternating} with stepsize $\eta = \tfrac{1}{5L}$ and negative momentum $\beta = -\tfrac{1}{2}$ satisfy
    \begin{align}
        \frac{1}{T}\sum_{t<T} \|\nabla f(z_t)\|^2 \leq 
        \frac{12\|z_0 - z^*\|^2}{\eta^2 T}.
    \end{align}
\end{theorem}
This result can be rewritten as $\mathbb{E}[\|\nabla f(z_\tau)\|^2] \lesssim \frac{1}{T}$ at a random stopping time $\tau$ chosen uniformly from $0,1,...,T-1$. Importantly, this does \emph{not} require taking an average of the iterates $z_t$, a key benefit for practical settings such as training GANs (see the discussions in e.g.,~\citep{mertikopoulos2019optimistic,hsieh2021limits} ).

The rate $\mathcal{O}(1/T)$ in~\cref{thm:cc-long} matches popular algorithms such as extragradient \citep{gorbunov2022extragradient,cai2022tight} and optimistic GDA \citep{Gorbunov_Taylor_Gidel, cai2022tight}. An accelerated rate of $\mathcal{O}(1/T^2)$ can be achieved by sophisticated algorithms \citep{yoon2021accelerated}, and it is an interesting question whether negative momentum can similarly lead to acceleration. By providing the first global convergence result for negative momentum,~\cref{thm:cc-long} opens the door to such refined questions. 

Our second result answers~\cref{q:fast}:

\begin{theorem}[Convergence for strongly-convex-strongly-concave objectives]\label{thm:scsc-long}
    Let $f$ be an $L$-smooth, $\mu$-strongly-convex-strongly-concave function with saddle point $z^* = (x^*,y^*)$. For any dimensions $d_x,d_y \in \N$, any initialization $z_0 = (x_0,y_0)$, and any number of iterations $T$, the iterates of alternating GDA~\eqref{eq:gda-momentum-alternating} with stepsize $\eta = \tfrac{1}{5L}$ and negative momentum $\beta = -\tfrac{1}{2}$ satisfy
    \begin{align}
        \|z_T - z^*\|^2 \leq 6( 1 - \eta \mu)^T \|z_0 - z^*\|^2.
    \end{align}
     In particular, 
     $\|z_T - z^*\|^2 \leq \eps \|z_0 - z^*\|^2$ after $T=\mathcal{O}(\kappa \logeps)$ iterations, where $\kappa := L/\mu$ denotes the condition number.
\end{theorem}

This rate $\mathcal{O}(\kappa \logeps)$ improves over the $\Omega(\kappa^{1.5}\logeps)$ lower bound for \emph{simultaneous} GDA with negative momentum~\citep{Zhang_Wang_2021}. Moreover, this rate is optimal among the class of first-order algorithms~\citep{Azizian_Mitliagkas_Lacoste-Julien_Gidel_2020} and matches classic algorithms such as extragradient and optimistic GDA~\citep{Mokhtari_Ozdaglar_Pattathil_2020}.

\begin{figure}
    \centering 
    \includegraphics[width=0.6\linewidth]{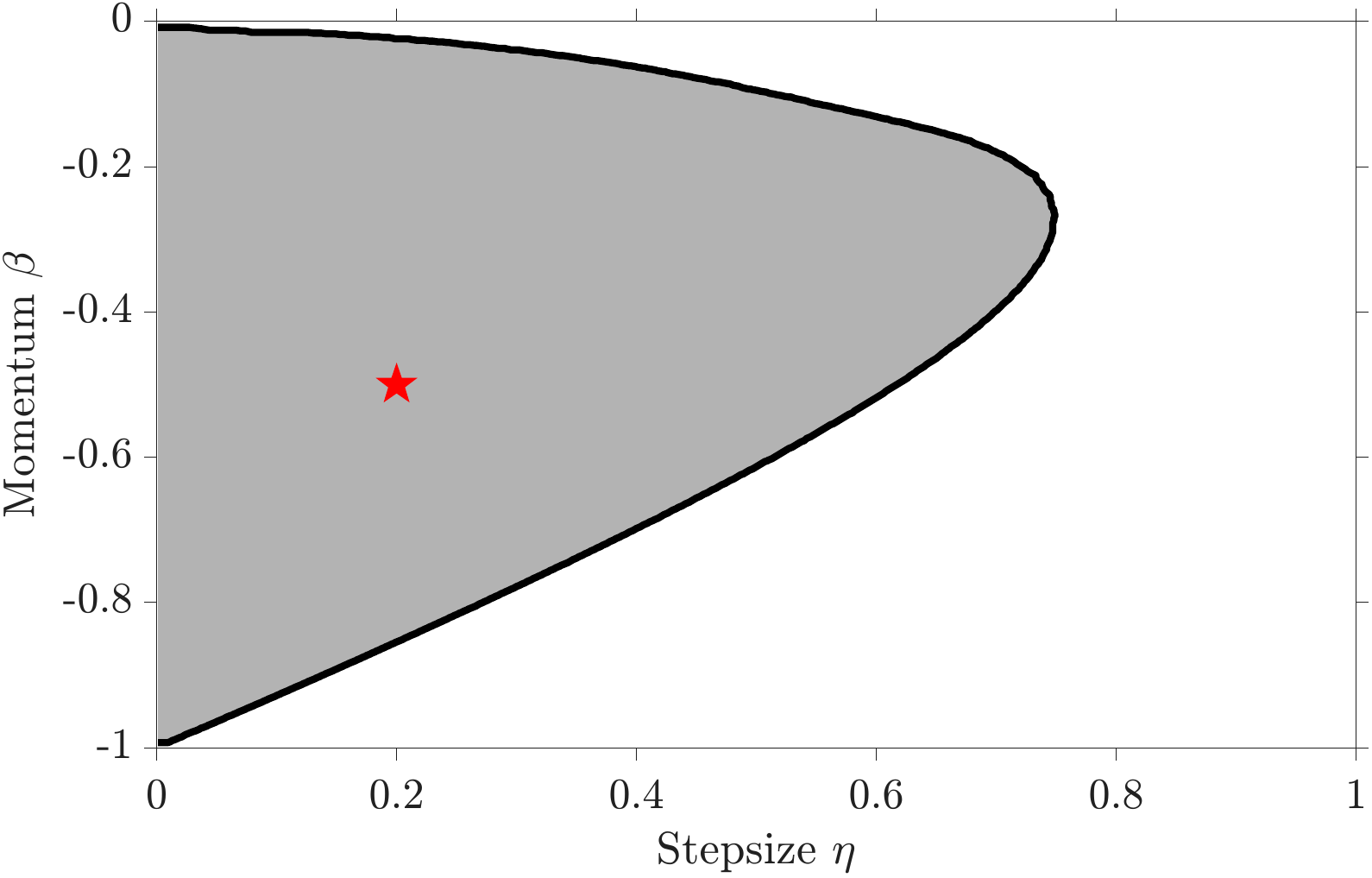}
   \caption{\footnotesize 
    The shaded region represents the set of parameters $(\eta,\beta)$ for which it is possible to certify convergence on 1-smooth, convex-concave functions ($\mu=0$) via the semidefinite programming approach described in~\cref{sec:proof} (i.e., the parameters for which convergence can be certified by an identity of the form~\eqref{eq:identity}).
   The boundary's roughness is due to numerical instability of SDP solvers (note however that our actual proofs are fully symbolic and rigorous). The red star marks the specific choice $(\eta,\beta) = (1/5,-1/2)$ considered for simplicity in \cref{thm:cc-long,thm:scsc-long}.}
    \label{fig:parameters}
\end{figure}

\begin{remark}[Choice of parameters]
    For concreteness and simplicity of exposition, we choose explicit constants for the stepsize $ \eta$ and momentum $\beta$ in~\cref{thm:scsc-long,thm:cc-long}. No effort has been made to optimize constants. The algorithm converges for other parameter choices, see~\cref{fig:parameters}. 
\end{remark}

\section{Progress lemma and its implications}\label{sec:lemma}

We establish~\cref{thm:scsc-long,thm:cc-long} via a progress lemma which identifies a quadratic Lyapunov function $\xi_t^\top \bm Q \xi_t$ under which negative momentum makes significant progress in each iteration. The state $\xi_t \in \R^{3(d_x + d_y)}$ includes the position, gradient\footnote{Note that the state $\xi_t$ does not include $\nabla_y f(x_{t+1},y_t)$ which is used in the second part of the update~\eqref{eq:gda-momentum-alternating}. This is for simplicity as this concise state $\xi_t$ captures enough information to show the desired convergence.}, and momentum at time $t$ in both the $x$ and $y$ coordinates:
$$\xi_t := (x_t-x^*,\, \nabla_x f(x_t,y_t),\, v_t, \, y_t-y^*,\, \nabla_y f(x_t,y_t),\, w_t ).$$ 
Here the momentum $(v_t,w_t)$ plays the role of $(x_t - x_{t-1}, y_t - y_{t-1})$. It is initialized to $(v_0,w_0) = (0,0)$ due to the standard convention of initializing $x_{-1} :=x_0$, and $y_{-1}:= y_0$.

The matrix $\bm Q $ can be interpreted as a coordinate system which extracts quadratic statistics of the state $\xi_t$ relevant for establishing convergence. We take $\bm Q$ to be block-diagonal of the form $$\bm Q = \begin{bmatrix}  \bm Q_x \otimes \bm I_{d_x} & \bm 0\\ \bm 0 & \bm Q_y \otimes \bm I_{d_y} \end{bmatrix}\,,$$ where $\bm Q_x, \bm Q_y \in \R^{3 \times 3}$ are positive definite.

For simplicity, we state this progress lemma for $1$-smooth functions and fixed parameters $\eta = 1/5$ and $\beta = -1/2$. This extends trivially to $L$-smooth functions via rescaling, see~\cref{rem:rescale} below.

\begin{lemma}[Progress lemma for negative momentum]\label{lem:progress}
    Let $\mu \in [0,1)$. Suppose $f$ is a $1$-smooth, $\mu$-strongly-convex-strongly-concave function with saddle point $z^* = (x^*,y^*)$. Consider the state vectors $\xi_t$ generated by running alternating GDA~\eqref{eq:gda-momentum-alternating} with stepsize $\eta = \frac{1}{5}$ and negative momentum $\beta = -\frac{1}{2}$. There exists $\bm Q$ such that for all iterations $t$, it holds that
    \begin{align} \label{eq:1step}
    (1-\mu/5)\xi_t^\top \bm Q \xi_t  - \xi_{t+1}^\top \bm Q \xi_{t+1} \geq   (1-\mu)\|\nabla f(z_t)\|^2
    \end{align}    
    and
    \begin{align}\label{eq:1step-Q}
        50 \|z_t -z^*\|^2 \leq \xi_t^\top \bm Q \xi_t \leq 150\|\xi_t\|^2 \,.
    \end{align}
\end{lemma}

This lemma is written in a convenient way so that it applies to both the strongly-convex-strongly-concave setting ($\mu > 0$) and the convex-concave setting ($\mu = 0$). The key inequality~\eqref{eq:1step} quantifies the decrease of the Lyapunov function $\xi_t^\top \bm Q \xi_t$, implying exponential convergence rates when $\mu >0$ (due to the multiplicative decrease from the contraction) and polynomial rates when $\mu=0$ (due to the additive decrease from the gradient term). The second inequality~\eqref{eq:1step-Q} ensures that the Lyapunov function can be related to quantities such as $\|z_t - z^*\|^2$, so that our final convergence results can be stated in terms of standard interpretable metrics.

\begin{remark}[Rescaling $L=1$ without loss of generality]\label{rem:rescale} Running the algorithm~\eqref{eq:gda-momentum-alternating} with momentum $\beta$ and stepsize $\eta = 1/(5L)$ on an $L$-smooth, $\mu$-strongly-convex-strongly-concave function $f$ generates exactly the same trajectory as running it with stepsize $1/5$ on the rescaled function $\tilde f := f/L$,
which is $1$-smooth and $\tilde{\mu}$-strongly-convex-strongly-concave for $\tilde{\mu} := \mu/L = 1/\kappa$. This is why we state~\cref{lem:progress} for $L=1$ without loss of generality. 
\end{remark}

\subsection{Proof of~\cref{thm:cc-long} using the progress lemma}\label{ssec:main-cc}

By~\cref{rem:rescale}, it is equivalent to analyze the trajectory on the rescaled function $\tilde f = f/L$ using rescaled stepsize $\tilde \eta = 1/5$. 
Note that $\xi_t$ refers to the state of this rescaled trajectory, and thus includes information about $\nabla \tilde f(x_t,y_t)$ rather than $\nabla f(x_t,y_t)$.

Setting $\mu=0$ and telescoping the key inequality~\eqref{eq:1step} over $T$ iterations gives:
$$ \sum_{t<T} \| \nabla \tilde f(z_t)\|^2 \leq \xi_0^\top \bm Q \xi_0 - \xi_{T}^\top \bm Q \xi_{T} \,.$$
We drop the term $\xi_{T}^\top \bm Q \xi_{T} \geq 0$ and upper bound $\xi_0^\top \bm Q \xi_0 \leq 150 \|\xi_0\|^2 \leq 300 \|z_0-z^*\|^2$ using property~\eqref{eq:1step-Q} and then the fact that $\|\xi_0 \|^2 = \|z_0-z^*\|^2+\|\nabla \tilde f(z_0)\|^2  \leq 2\|z_0-z^*\|^2$ (by $1$-smoothness of $f$ and the standard initialization $z_{-1} = z_0$ which implies $v_0,w_0=0$). Plugging in $\nabla \tilde f = \nabla f/L$ gives  
$$ \frac{1}{L^2}\sum_{t<T} \| \nabla  f(z_t)\|^2 \leq 300 \|z_0-z^*\|^2\,.$$
Multiplying both sides by $\tfrac{L^2}{T}$ and recalling the stepsize choice $\eta = \tfrac{1}{5L}$ completes the proof.

\subsection{Proof of~\cref{thm:scsc-long} using the progress lemma}\label{ssec:main-scsc}

Again by~\cref{rem:rescale} we analyze the trajectory on the rescaled function $\tilde f := f/L$. The proof of~\cref{thm:scsc-long} then follows immediately from~\cref{lem:progress}:
$$ \|z_T-z^*\|^2 \leq \frac{1}{50} \xi_{T}^\top \bm Q \xi_{T} \leq \frac{1}{50}(1-\eta\mu)^T\xi_0^\top \bm Q \xi_0 \leq 3(1-\eta\mu)^T \|\xi_0\|^2 \leq 6(1-\eta\mu)^T \|z_0-z^*\|^2.$$
Above, the first and third steps are by the property~\eqref{eq:1step-Q} of $\bm Q$. The second step is by the key inequality~\eqref{eq:1step} in~\cref{lem:progress}, which implies the $1$-step contraction $\xi_{t+1}^\top \bm Q \xi_{t+1} \leq
(1 - \tilde{\mu}/5)
\xi_t^\top \bm Q \xi_t = (1 - \eta \mu) \xi_t^\top \bm Q \xi_t $ after dropping the gradient term. The final step is because $\|\xi_0 \|^2 \leq 2\|z_0-z^*\|^2$ as argued already in~\cref{ssec:main-cc}.

%% file: sections/proof_key_lemma.tex
\section{Proof of progress lemma}\label{sec:proof}

In the previous section, we showed that \cref{lem:progress} implies our main results; in this section, we prove this key progress lemma.

We prove~\cref{lem:progress} by establishing an identity of the form
\begin{equation}\label{eq:identity}
(1-\mu/5)\xi_t^\top \bm Q \xi_t  - \xi_{t+1}^\top \bm Q \xi_{t+1} -(1-\mu)\|\nabla f(z_t)\|^2 = \sum_{\alpha} \lambda_\alpha M_\alpha + S
\end{equation}
for an explicit matrix $\bm Q$ satisfying~\eqref{eq:1step-Q}. 
The left-hand side of this identity is the quantity that we seek to show is non-negative in order to prove the progress inequality~\eqref{eq:1step} in~\cref{lem:progress}.
The right-hand side consists of two types of non-negative quantities. 
The first term is the sum of non-negative multipliers $\lambda_\alpha \geq 0$ (that we construct explicitly) multiplied by certain explicit polynomial quantities $M_\alpha \geq 0$ of the iterates (that are non-negative by the smoothness and convexity-concavity 
properties of $f$). The second term $S$ is a sum-of-squares quadratic polynomial.

\subsection{Starting point and technical overview} 

Our starting point for proving~\eqref{eq:identity} is a general numerical approach, based on semidefinite programming (SDP), for searching for the quantities in identities of this form~\citep{taylor2018lyapunov,upadhyaya2025automated, Lessard_Recht_Packard_2016}. Intuitively, this search is an SDP feasibility problem because there are linear equality constraints and positive semidefinite (PSD) constraints in the relevant variables $\{\lambda_{\alpha}\}$, $\bm Q$, and $S$. Indeed, the linear equality constraints arise because both sides of the identity~\eqref{eq:identity} are quadratic polynomials (in the state space) whose coefficients depend linearly on the decision variables $\{\lambda_{\alpha}\}$, $\bm Q$, $S$. The PSD constraints arise because $S$ being a sum-of-squares quadratic is equivalent to it being a PSD quadratic form in the state space, and because $\bm Q$ satisfying the two-sided inequality~\eqref{eq:1step-Q} is equivalent to two PSD constraints.

However, while this connection to SDP provides a helpful starting point, there are three overarching challenges that we must address for our problem:

\paragraph{Challenge 1: choice of valid inequalities.} Establishing the identity~\eqref{eq:identity} requires first specifying the non-negative quantities $M_{\alpha} \geq 0$. This amounts to choosing which ``valid inequalities'' or ``proof system'' one uses to prove a convergence rate, or equivalently which structural properties of the objective $f$ to exploit. The most popular approach for choosing $M_{\alpha}$ in min-max optimization settings is based on the classical fact that if $f(x,y)$ is convex-concave, then $(\nabla_x f, -\nabla_y f)$ is a monotone operator; this monotonicity suggests simple quantities $M_{\alpha} \geq 0$. See for example \citep{Ryu_Taylor_Bergeling_Giselsson_2020}. However, a fundamental issue is that those quantities $M_{\alpha}$ are unavoidably oblivious to the asymmetry between $x$ and $y$ intrinsic to convex-concave functions $f(x,y)$. Without this structure, the proof system is too weak to establish the identity~\eqref{eq:identity}. We overcome this by using a more powerful proof system $\{M_{\alpha}\}$ that directly stems from the definition of convex-concavity. Details in \cref{ssec:lemma-multipliers}.

\paragraph*{Challenge 2: parametric solution.} The SDP approach numerically certifies the identity~\eqref{eq:identity}, but only for a fixed numerical choice of algorithm parameters (stepsize $\eta$ and momentum $\beta$) and  problem parameters (strong convexity $\mu$ and smoothness $L$). However, this is insufficient due to two important issues. First, it is a priori unclear how to find good algorithm parameters because the search for the stepsize $\eta$ and momentum $\beta$ cannot be tractably included in the SDP---indeed this search is in general non-convex, a notorious challenge in the PEP literature more broadly~\citep{taylor2024towards}. Second, \cref{lem:progress} requires a (symbolic) proof that applies to any value of the problem parameter $\mu$ (one can eliminate $L$ via rescaling, see~\cref{rem:rescale}). We overcome this by identifying convenient choices of $\eta$ and $\beta$ that enable SDP solutions and convergence proofs that are simple in two key ways. First, they require only a sparse subset of the multipliers $\{\lambda_{\alpha}\}$, in fact with only $3$ distinct values out of potentially $72$, which makes it significantly more tractable to identify symbolic solutions. Second, this enables choices of the multipliers $\lambda_\alpha$ and potential function $\bm Q$ that are \emph{independent} of the problem parameter $\mu \geq 0$, making the symbolic verification (described next) significantly simpler. Details in~\cref{ssec:lemma-construction}.

\paragraph*{Challenge 3: rigorous proof.} SDP solvers provide numerical outputs that may satisfy the feasibility constraints only approximately. However, proving~\cref{lem:progress} requires rigorously certifying both the linear coefficient-matching constraints defining~\eqref{eq:identity} as well as the PSD constraints related to~\eqref{eq:1step-Q} and the sum-of-squares property for $S$. 
Rounding numerical SDP solutions to rigorous proofs is a well-documented challenge in the PEP literature, see e.g., \citep{taylor2024towards}. We overcome this by leveraging techniques from the symbolic computational algebra community, for example we use Descartes' rule of signs and Sturm's method in order to certify PSD constraints (the most difficult of these issues for PEP). Details in~\cref{ssec:lemma-sos}.

\paragraph*{Organization of the remainder of the section.} Below we detail how we address these challenges and prove~\cref{lem:progress}. This analysis is inspired by the numerical output of the aforementioned SDP, but we emphasize that our proof is fully rigorous (not numerical) and can be read entirely by itself (the purpose of the above discussion is primarily to provide insight for how this proof was obtained). Specifically, below in~\cref{ssec:lemma-multipliers} we specify the choice of valid inequalities $\{M_{\alpha}\}$; in~\cref{ssec:lemma-construction} we explicitly construct all other quantities $\{\lambda_{\alpha}\}$, $\bm Q$, $S$ in the progress identity~\eqref{eq:identity}; in~\cref{ssec:lemma-sos} we prove that $S$ satisfies the desired sum-of-squares property; and finally in~\cref{ssec:lemma-proof} we combine these ingredients to prove~\cref{lem:progress}.

\subsection{Valid inequalities}\label{ssec:lemma-multipliers}

Here we specify the non-negative quantities $\{M_{\alpha}\}$ that we use in the progress identity~\eqref{eq:identity}. As explained above, this choice of $\{M_{\alpha}\}$ can be interpreted as the choice of which structural properties of $f$ (a.k.a.\ which ``valid inequalities'') we exploit in our convergence analysis. We use three types of quantities, defined for $i,j,k,l \in \{t,t+1,*\}$:

\begin{itemize}
    \item \textbf{Smoothness.} Recall that $f$ being $1$-smooth amounts to $\nabla f$ being $1$-Lipschitz (see~\cref{def:smooth}). Thus the following quantities are clearly non-negative:
    $$\Msmooth(x_i,y_j, x_k,y_l) := \|(x_i,y_j)-(x_k,y_l)\|^2-\|\nabla f(x_i,y_j)-\nabla f(x_k,y_l)\|^2.$$
    \item \textbf{Convexity.} Since $f(\cdot,\cdot)$ is $\mu$-strongly-convex-strongly-concave and $1$-smooth, the restricted function $f(\cdot,y_k)$ is $\mu$-strongly convex and $1$-smooth for any $y_k$. By \cref{lem:cocoercivity}, this implies non-negativity of the associated co-coercivities:
  \begin{align*}
    &\Mconvex(x_i,x_j, y_k)
    :=
    C_{f(\cdot, y_k)}(x_i,x_j)\,.
\end{align*}
    \item \textbf{Concavity.} Analogously, the restricted function $f(x_k,\cdot)$ is $\mu$-strongly-concave and $1$-smooth, hence the co-coercivities for $-f(x_k,\cdot)$ are non-negative:
    \begin{align*}
        \Mconcave(y_i,y_j,x_k) := C_{-f(x_k,\cdot)}(y_i,y_j)\,.
    \end{align*}
\end{itemize}

Note that all three quantities are  polynomials which are linear (or zero) in the function values and quadratic in the iterates and gradients. Indeed, by expanding the definition of the co-coercivities (\cref{lem:cocoercivity}), the latter two quantities can be written explicitly as 
\begin{align*}
    &\Mconvex(x_i,x_j, y_k) = f(x_i,y_k)-f(x_j,y_k)-\nabla_x f(x_j,y_k)^\top(x_i-x_j)
    \\ &\qquad\qquad -\frac{\mu}{2}\|x_i-x_j\|^2 - \frac{1}{2(1-\mu)}\|\nabla_x f(x_i,y_k) - \nabla_x f(x_j,y_k) - \mu(x_i-x_j)\|^2\,, 
    \\
    &\Mconcave (y_i,y_j,x_k) =   f(x_k,y_j)-f(x_k,y_i) +\nabla_y f(x_k,y_j)^\top(y_i-y_j)
    \\ &\qquad\qquad -\frac{\mu}{2}\|y_i-y_j\|^2 - \frac{1}{2(1-\mu)}\|\nabla_y f(x_k,y_j) - \nabla_y f(x_k,y_i) - \mu(y_i-y_j)\|^2\,.
\end{align*}

These non-negative quantities $M_{\alpha}$ are also used to prove convergence rates in \citep{shugart_negative_2025, Lee_Cho_Yun_2024}, albeit for different algorithms and purposes. We highlight two crucial features of this choice:

\begin{enumerate}
    \item \textbf{Convexity-concavity rather than monotonicity.} More standard is not to use  $\Mconvex$ and $\Mconcave$, instead to use only the simpler $\Mmonotone(z_i,z_j) := \langle F(z_i) - F(z_j), z_i - z_j \rangle$ where $z_i := (x_i,y_i)$, $z_j := (x_j,y_j)$, and $F := (\nabla_x f, -\nabla_y f)$ is the saddle-point operator associated to $f$. The motivation is that convexity-concavity of $f$ implies monotonicity of the operator $F$, which ensures non-negativity of $\Mmonotone$. See for example~\citep{Ryu_Taylor_Bergeling_Giselsson_2020, Zhang_Grosse}. However, convexity-concavity of $f$ is a strictly stronger property than monotonicity of $F$ when enforced only at a finite collection of points (as done here), which is why $\Mconvex$ and $\Mconcave$ are strictly more powerful inequalities than $\Mmonotone$ and enable proving much faster convergence rates. 

    \item \textbf{Gridded inequalities.} More standard is to use valid inequalities which depend on the iterates only through $(x_i,y_i)$ for $i \in \{t,t+1,*\}$. In contrast,  we consider valid inequalities of the above three types for  $i,j,k,l\in \{t,t+1,*\}$. This amounts to enforcing these inequalities on a ``grid'' of (concatenated) points, many of which are not visited by the algorithm yet are helpful to proving convergence analyses. For example, $\Mconcave(y_t,y_{t+1},x_t)$ includes information about $f$ at $(x_t,y_{t+1})$, and $\Mconcave(y^*,y_t,x_{t+1})$ includes information about $f$ at $(x_{t+1},y^*)$. This gridding enables exploiting global structural properties of $f$ that the aforementioned standard approach cannot.
\end{enumerate}
Both features enable proving stronger rates, and the combination is needed to prove~\cref{lem:progress}.

\subsection{Construction of quantities in progress identity~\eqref{eq:identity}}\label{ssec:lemma-construction}

Here we explicitly construct the quantities $\bm Q$, $\bm S$, and $\{\lambda_{\alpha}\}$ that define the progress identity~\eqref{eq:identity}.

\par Letting $\otimes$ denote the Kronecker product, we define $\bm Q$ as
\begin{equation}
    \bm Q := \begin{bmatrix}  \bm Q_x \otimes \bm I_{d_x} & \bm 0\\ \bm 0 & \bm Q_y \otimes \bm I_{d_y} \end{bmatrix}, \quad \bm Q_x := \begin{bmatrix}
        120 & -\frac{81}{10} & -40\\ -\frac{81}{10} &1&4\\-40&4&40
    \end{bmatrix}, \quad \bm Q_y := \begin{bmatrix}
        120 & 0 & -40\\ 0 &1&-4\\-40&-4&40
    \end{bmatrix}.
\end{equation}
Note that $\bm Q$ satisfies~\eqref{eq:1step-Q} since a straightforward eigenvalue calculation shows that $50 \bm E_{11} \preceq \bm Q_x, \bm Q_y \preceq 150 \bm I$, where $\bm E_{11}$ denotes the $3 \times 3$ matrix with first entry $1$ and all other entries $0$.

Next, we define the multipliers $\{\lambda_{\alpha}\}$. Although the proof system in~\cref{ssec:lemma-multipliers} allows for  arbitrary non-negative combinations $\sum_{\alpha} \lambda_\alpha M_\alpha$ of all $72$ valid inequalities $M_{\alpha}$, an appealing aspect of our analysis is that we use only a small subset of these quantities to prove the desired identity.
Explicitly, letting $\lambda_1 = 4, \lambda_2 = \frac{79}{5}, \lambda_3 = \frac{81}{5}$, the non-negative combination of valid inequalities we use is
\begin{align*}
    \sum_{\alpha} \lambda_\alpha M_\alpha &=  \lambda_1 \Msmooth(x_t,y_t,x_{t+1},y_{t+1})\\  &\ + \lambda_2\left( \Mconvex(x_t,x_{t+1},y_{t}) + \Mconvex(x^*,x_{t},y_{t}) + \Mconcave(y_t,y^*,x^*)\right)\\ &\ + \lambda_3\left( \Mconvex(x^*,x_{t+1},y_{t+1})+ \Mconcave(y_{t+1},y^*,x^*)+ \Mconcave(y_t,y_{t+1},x_t)  \right)\\ &\ + (\lambda_2+\lambda_3) \left(\Mconvex(x_{t+1},x^*,y^*)+ \Mconcave(y^*,y_t,x_{t+1})\right)
\end{align*}

Finally, we define $S$ as the residual in~\eqref{eq:identity} so that the identity holds by construction:
\begin{align}\label{eq:def-S}
    S := (1-\mu/5)\xi_t^\top \bm Q \xi_t  - \xi_{t+1}^\top \bm Q \xi_{t+1} -(1-\mu)\|\nabla f(z_t)\|^2 - \sum_{\alpha} \lambda_\alpha M_\alpha\,.
\end{align}
It remains to  prove that this $S$ satisfies the desired property; we do this next.

\subsection{Proof that $S$ is sum-of-squares}\label{ssec:lemma-sos}

As described above, in order to use the progress identity~\eqref{eq:identity} to prove~\cref{lem:progress}, we need to show that $S$ is non-negative. We accomplish this by showing that $S$ is a sum of squares for any $\mu \in [0,1)$. We begin by simplifying the definition of $S$ in~\eqref{eq:def-S}.
An algebraically tedious but conceptually straightforward calculation shows that $S$ is a quadratic form which can be expressed as
\begin{align}\label{eq:S-explicit}
S = \frac{1}{2(1-\mu)}\Xi_t^\top \begin{bmatrix} \bm S_x \otimes \bm I_{d_x} & \bm 0 \\ \bm 0 & \bm S_y \otimes \bm I_{d_y}  \end{bmatrix} \Xi_t\,.
\end{align}
The matrices $\bm S_x, \bm S_y$ are complicated but fully explicit (see Appendix~\ref{app:sos} for the closed-form) and have entries which depend quadratically in $\mu$.
The vector $\Xi_t$ is an extended state space that incorporates information about $\nabla f$ at certain additional points that are not iterates of the algorithm:  
\begin{multline*}
\Xi_t := \left( x_t-x^*, \nabla_x f(x^*,y_t), \nabla_x f(x^*,y_{t+1}), \nabla_x f(x_{t+1},y^*), \nabla_x f(x_t,y_t), \nabla_x f(x_{t+1},y_t), \nabla_x f(x_{t+1},y_{t+1}), v_t,\right.\\ \left.
y_t-y^*,  \nabla_y f(x^*,y_t), \nabla_y f(x^*,y_{t+1}), \nabla_y f(x_{t+1},y^*), \nabla_y f(x_t,y_t), \nabla_y f(x_{t+1},y_t), \nabla_y f(x_{t+1},y_{t+1}) ,w_t \right).
\end{multline*}

By virtue of the expression~\eqref{eq:S-explicit}, the following lemma suffices to show that $(1-\mu)S$ is a sum-of-squares polynomial (and therefore $S$ is non-negative) for any fixed value of $\mu \in [0,1)$.

\begin{lemma}\label{lem:sos}
    For any $\mu \in [0,1)$, it holds that $\bm S_x, \bm S_y \succeq 0$.
\end{lemma}

For brevity, we sketch the main idea of the proof here and defer full details to Appendix~\ref{app:sos}. We prove~\cref{lem:sos} by showing that for all $\mu \in [0,1)$, the characteristic polynomials of $\bm S_x, \bm S_y$ have only non-negative roots. We establish this via Descartes' rule of signs: it suffices to show that these characteristic polynomials have coefficients which alternate in sign. Each coefficient is itself a polynomial in $\mu$, so we must verify the alternating sign pattern holds for all $\mu \in [0,1)$. This can be rigorously proved by first checking this pattern for $\mu = 0$ and then confirming, using standard symbolic algebra techniques such as Sturm's method, that the relevant coefficient polynomials have no roots in $(0,1)$ and thus never change sign.

\subsection{Concluding the proof of~\cref{lem:progress}}\label{ssec:lemma-proof}

Consider $\{M_{\alpha}\}$ as defined in~\cref{ssec:lemma-multipliers}, and consider $\{\lambda_{\alpha}\}$, $\bm Q$, and $S$ as defined in~\cref{ssec:lemma-construction}. By construction of $S$, the identity~\eqref{eq:identity} holds. Next, note that $\lambda_{\alpha}$, $M_{\alpha}$, and $S$ are all non-negative (the former two by construction, and the latter by~\cref{lem:sos}). Thus the right-hand side of~\eqref{eq:identity} is non-negative, which implies the desired progress bound~\eqref{eq:1step}. Finally, the property~\eqref{eq:1step-Q} of $\bm Q$ is immediate from the construction of $\bm Q$ in~\cref{ssec:lemma-construction}.

%% file: sections/appommitteddetails.tex
\section{Deferred details from~\cref{ssec:lemma-sos}}\label{app:sos}

Here we provide the remaining details for~\cref{ssec:lemma-sos} and in particular prove~\cref{lem:sos}. See~\cref{ssec:lemma-sos} for an overview of the proof strategy.

\subsection{Explicit expressions}\label{app:sos-explicit}

We begin by explicitly stating the coefficient matrices $\bm S_x$ and $\bm S_y$ in the factorization~\eqref{eq:S-explicit} of $S$:

\small 

$$ \bm S_x = \begin{bmatrix}
16\mu + 48\mu^2 & \frac{79}{5}\mu & \frac{81}{5}\mu & -32\mu & -\frac{111}{5}\mu - \frac{81}{25}\mu^2 & 0 & -\frac{81}{5}\mu & -\frac{81}{10}\mu - 16\mu^2 \\
\frac{79}{5}\mu & \frac{79}{5} & 0 & 0 & -\frac{79}{5} & 0 & 0 & 0 \\
\frac{81}{5}\mu & 0 & \frac{81}{5} & 0 & -\frac{81}{25}\mu & 0 & -\frac{81}{5} & -\frac{81}{10}\mu \\
-32\mu & 0 & 0 & 32 & \frac{32}{5}\mu & 0 & 0 & 16\mu \\
-\frac{111}{5}\mu - \frac{81}{25}\mu^2 & -\frac{79}{5} & -\frac{81}{25}\mu & \frac{32}{5}\mu & \frac{822}{25} - \frac{86}{25}\mu - \frac{8}{5}\mu^2 & -\frac{316}{25} & -\frac{32}{5} + \frac{241}{25}\mu & -\frac{44}{5} + \frac{57}{10}\mu + \frac{8}{5}\mu^2 \\
0 & 0 & 0 & 0 & -\frac{316}{25} & \frac{79}{5} & 0 & \frac{79}{10} \\
-\frac{81}{5}\mu & 0 & -\frac{81}{5} & 0 & -\frac{32}{5} + \frac{241}{25}\mu & 0 & \frac{111}{5} - 6\mu & 4 + \frac{41}{10}\mu \\
-\frac{81}{10}\mu - 16\mu^2 & 0 & -\frac{81}{10}\mu & 16\mu & -\frac{44}{5} + \frac{57}{10}\mu + \frac{8}{5}\mu^2 & \frac{79}{10} & 4 + \frac{41}{10}\mu & 38 - 38\mu + 16\mu^2
\end{bmatrix}$$

$$\bm S_y = \begin{bmatrix}
16\mu + 48\mu^2 & \frac{79}{5}\mu & \frac{81}{5}\mu & -32\mu & 0 & \frac{881}{25}\mu & 0 & \frac{79}{10}\mu - 16\mu^2 \\
\frac{79}{5}\mu & \frac{79}{5} & 0 & 0 & 0 & 0 & 0 & 0 \\
\frac{81}{5}\mu & 0 & \frac{81}{5} & 0 & 0 & \frac{81}{25}\mu & 0 & -\frac{81}{10}\mu \\
-32\mu & 0 & 0 & 32 & 0 & -32 & 0 & 0 \\
0 & 0 & 0 & 0 & 8 - \frac{32}{5}\mu - \frac{8}{5}\mu^2 & 0 & -8 + 8\mu & -8 + \frac{48}{5}\mu - \frac{8}{5}\mu^2 \\
\frac{881}{25}\mu & 0 & \frac{81}{25}\mu & -32 & 0 & \frac{1037}{25} + \frac{192}{125}\mu & -\frac{284}{25} - \frac{8}{5}\mu & \frac{84}{5} - \frac{597}{50}\mu \\
0 & 0 & 0 & 0 & -8 + 8\mu & -\frac{284}{25} - \frac{8}{5}\mu & \frac{111}{5} - 6\mu & -\frac{121}{10} + 4\mu \\
\frac{79}{10}\mu - 16\mu^2 & 0 & -\frac{81}{10}\mu & 0 & -8 + \frac{48}{5}\mu - \frac{8}{5}\mu^2 & \frac{84}{5} - \frac{597}{50}\mu & -\frac{121}{10} + 4\mu & 38 - \frac{459}{10}\mu + 16\mu^2
\end{bmatrix}$$

\normalsize

In the proof of~\cref{lem:sos} below, we make use of the characteristic polynomials $p_x(\zeta;\mu) := \det\left( \bm S_x -\zeta \bm I\right)$ and $p_y(\zeta;\mu) := \det\left( \bm S_y -\zeta \bm I\right)$ of these matrices. Since each entry of $\bm S_x$ and $\bm S_y$ is a polynomial in $\mu$, these characteristic polynomials are bivariate in $\zeta,\mu$ and can be written as
\begin{align}\label{eq:characteristic}
    p_x(\zeta;\mu) = \sum_{j=0}^8 c_{x,j}(\mu) \zeta^j,\qquad p_y(\zeta;\mu) = \sum_{j=0}^8 c_{y,j}(\mu) \zeta^j,
\end{align}
where $c_{x,j}(\mu)$ and $c_{y,j}(\mu)$ are polynomials of $\mu$. These coefficients can be computed explicitly via an algebraically tedious but conceptually simple calculation. We state these coefficients below in a concatenated (matrix) form to make it easier to visually see the sign changes, which is the key way in which we use these coefficients below.
For $i,j \geq 0$, the coefficient of $\mu^i$ in $c_{x,j}(\mu)$ is the $(i+1,j+1)$-th entry of $\bm C_x$. Similarly for $\bm C_y$.

\begin{equation*}\label{eq:px}
\bm C_x=\begin{bsmallmatrix}
0 & -\frac{3016949328}{3125} & \frac{3925293412889}{31250} & -\frac{315920718241}{6250} & \frac{410082261313}{62500} & -\frac{4859926311}{12500} & \frac{29183151}{2500} & -\frac{4322}{25} & 1 \\
\frac{48271189248}{3125} & -\frac{32554037377816}{15625} & \frac{19024674172661}{31250} & -\frac{116858962492}{3125} & -\frac{4944502477}{6250} & \frac{78866357}{625} & -\frac{444266}{125} & \frac{786}{25} & 0 \\
\frac{118119041994672}{390625} & \frac{218825704369549}{31250} & -\frac{1628209078050351}{781250} & \frac{24840599510927}{156250} & -\frac{20682996319}{12500} & -\frac{647480389}{3125} & \frac{4380333}{625} & -\frac{312}{5} & 0 \\
-\frac{570636561473136}{390625} & -\frac{6749359461879353}{781250} & \frac{794651282514111}{390625} & -\frac{14880446763889}{156250} & -\frac{36988061198}{15625} & \frac{48458491}{250} & -\frac{318582}{125} & 0 & 0 \\
\frac{980677852825824}{390625} & \frac{1863197504635273}{390625} & -\frac{144439831610443}{156250} & \frac{5736968309049}{156250} & \frac{26987184513}{31250} & -\frac{315937633}{6250} & \frac{247839}{625} & 0 & 0 \\
-\frac{789913089425376}{390625} & -\frac{95808733158877}{78125} & \frac{226228514254623}{781250} & -\frac{1259177102927}{78125} & \frac{2732915364}{15625} & \frac{1390428}{625} & 0 & 0 & 0 \\
\frac{293540264356464}{390625} & \frac{152705311358489}{781250} & -\frac{42651897368459}{781250} & \frac{609107362377}{156250} & -\frac{3267679759}{31250} & \frac{590096}{625} & 0 & 0 & 0 \\
-\frac{37821406934448}{390625} & -\frac{13840030176837}{781250} & \frac{377986454943}{78125} & -\frac{9381626397}{31250} & \frac{3540576}{625} & 0 & 0 & 0 & 0
\end{bsmallmatrix}
\end{equation*}
\begin{equation*}\label{eq:py}
\bm C_y=\begin{bsmallmatrix}
0 & -\frac{1288809792}{15625} & \frac{60597876882}{15625} & -\frac{246249337657}{12500} & \frac{246345620709}{62500} & -\frac{754682859}{2500} & \frac{26859791}{2500} & -\frac{4342}{25} & 1 \\
\frac{20620956672}{15625} & -\frac{5264240621856}{78125} & \frac{28547368182704}{78125} & -\frac{6613389982081}{156250} & \frac{4107920378}{15625} & \frac{15227937}{125} & -\frac{2794858}{625} & \frac{10191}{250} & 0 \\
\frac{6278652931968}{78125} & -\frac{240098071662172}{390625} & -\frac{652932152553913}{781250} & \frac{32815526836777}{312500} & -\frac{29665687223}{15625} & -\frac{1048066121}{6250} & \frac{16482201}{2500} & -\frac{312}{5} & 0 \\
-\frac{59961927583872}{390625} & \frac{1112181291312788}{390625} & \frac{449869267325923}{781250} & -\frac{9937608117483}{156250} & -\frac{10725965407}{15625} & \frac{450804572}{3125} & -\frac{1439996}{625} & 0 & 0 \\
-\frac{758510810673408}{1953125} & -\frac{8360973992985416}{1953125} & -\frac{2871838767839}{31250} & \frac{974565108823}{31250} & -\frac{8239681301}{62500} & -\frac{96310742}{3125} & \frac{10176}{25} & 0 & 0 \\
\frac{2827722585238272}{1953125} & \frac{5710819702531544}{1953125} & -\frac{9205535476859}{156250} & -\frac{1972122788511}{156250} & \frac{4290414752}{15625} & \frac{2966208}{3125} & 0 & 0 & 0 \\
-\frac{3367757488007808}{1953125} & -\frac{1801214735526316}{1953125} & \frac{3957829549436}{78125} & \frac{6639541174}{3125} & -\frac{1620729152}{15625} & \frac{23552}{25} & 0 & 0 & 0 \\
\frac{1850723490483072}{1953125} & \frac{235934676676004}{1953125} & -\frac{4122192150876}{390625} & -\frac{4024869248}{78125} & \frac{13142016}{3125} & 0 & 0 & 0 & 0 \\
-\frac{437038196779008}{1953125} & -\frac{14342187792576}{1953125} & \frac{547701173056}{390625} & -\frac{327867392}{15625} & 0 & 0 & 0 & 0 & 0 \\
\frac{5025222955008}{390625} & \frac{225511122816}{390625} & -\frac{927148032}{15625} & 0 & 0 & 0 & 0 & 0 & 0
\end{bsmallmatrix}
\end{equation*}

\subsection{Proof of~\cref{lem:sos}}

    Let $p_x(\zeta;\mu) := \det\left( \bm S_x -\zeta \bm I\right)$ and $p_y(\zeta;\mu) := \det\left( \bm S_y -\zeta \bm I\right)$ denote the characteristic polynomials of $\bm S_x$ and $\bm S_y$, respectively. (See~\eqref{eq:characteristic} for their explicit expressions in terms of coefficients $c_{x,j}(\mu)$ and $c_{y,j}(\mu)$.)
    Since $\bm S_x$ and $\bm S_y$ are symmetric, their eigenvalues are real, hence all roots of $p_x(\cdot;\mu)$ and $p_y(\cdot;\mu)$ are real for all $\mu$. To prove the lemma, we show that these real roots are in fact non-negative. By Descartes' rule of signs\footnote{We recall here the relevant version of Descartes' rule of signs (see e.g.~\citep[Theorem 2.33]{basu2006algorithms}). Let $q(\zeta) = \sum_{k=0}^n c_k \zeta^k$ have real coefficients and real roots. Suppose that for some $m \ge 0$, the low-order coefficients $c_k = 0$ for all $k < m$, and the high-order coefficients $c_m, c_{m+1}, \dots, c_n $ are non-zero and alternate in sign. Then all roots of $q$ are non-negative.}, it suffices to show that for each $\mu\in [0,1)$, the non-zero coefficients of $p_x(\zeta;\mu)$ and $p_y(\zeta;\mu)$ alternate in sign.

    First, for $\mu = 0$, Descartes' rule of signs holds by inspection of $c_{x,j}(0)$ and $c_{y,j}(0)$. See the first row of $\bm C_x$ and $\bm C_y$ in~\cref{app:sos-explicit} above.

    Next, to apply Descartes' rule of signs for $\mu > 0$, we first deal with a slight nuance: the coefficients $c_{x,0}$ and $c_{y,0}$ vanish at $\mu=0$ (see the top-left entries of $\bm C_x$ and $\bm C_y$). To ensure that the alternating sign pattern persists for small positive $\mu$, observe that in both $c_{x,0}$ and $c_{y,0}$, the linear coefficient of $\mu$ is positive (see the (2,1)-th entries of $\bm C_x$ and $\bm C_y$), which in particular is the opposite of the sign of the constant term in $c_{x,1}$ and $c_{y,1}$ (see the (1,2)-th entries). All other coefficients are already nonzero at $\mu=0$, hence by continuity their signs remain unchanged for sufficiently small $\mu>0$, and thus the alternating sign pattern persists for sufficiently small $\mu > 0$.

    Now we are ready to apply Descartes' rule of signs simultaneously for all $\mu \in (0,1)$. To do this, it suffices to show that $c_{x,j}(\mu)$ and $c_{y,j}(\mu)$ do not change sign on the interval $\mu \in (0,1)$. That is, it suffices to show that these coefficient polynomials $c_{x,j}(\mu)$ and $c_{y,j}(\mu)$ have no roots on the interval.  This can be rigorously proven using standard techniques from symbolic computer algebra such as Sturm's Theorem 
    (see e.g. \citep[Theorem 2.50]{basu2006algorithms}). For brevity of exposition and the convenience of the reader, we provide a short Mathematica script that validates this in a  rigorous symbolic manner~\citep{MathematicaURL}.